\documentclass[a4paper, 12pt]{article}
\usepackage[utf8]{inputenc}
\usepackage{amsmath}
\usepackage{amssymb}
\usepackage{multirow}
\usepackage{fullpage}
\usepackage{graphicx}
\usepackage{amsthm}
\numberwithin{equation}{section}
\DeclareMathOperator{\ed}{d}
\DeclareMathOperator{\vect}{vect}
\DeclareMathOperator{\diam}{diam}

\newcommand{\wm}{\mathcal{W}}
\newcommand{\drm}{\mathcal{C}}

\newcommand{\im}{\mathfrak{I}}
\newtheorem{proposition}{Proposition}[section]

\newtheorem{theorem}[proposition]{Theorem}
\newtheorem{lemma}[proposition]{Lemma}
\theoremstyle{definition}
\newtheorem{definition}[proposition]{Definition}

\theoremstyle{remark}
\newtheorem{remark}[proposition]{Remark}

\title{New degrees of freedom for differential forms on cubical meshes}
\date{}
\author{Jonni Lohi\\Faculty of Information Technology\\University of Jyväskylä}

\begin{document}

\maketitle

\begin{abstract}
We consider new degrees of freedom for higher order differential forms on cubical meshes. The approach is inspired by the idea of Rapetti and Bossavit to define higher order Whitney forms and their degrees of freedom using small simplices. We show that higher order differential forms on cubical meshes can be defined analogously using small cubes and prove that these small cubes yield unisolvent degrees of freedom. Significantly, this approach is compatible with discrete exterior calculus and expands the framework to cover higher order methods on cubical meshes, complementing the earlier strategy based on simplices.
\end{abstract}

\section{Introduction}

Finite element exterior calculus \cite{arnold2006} highlights the importance of suitable finite element spaces in discretisations of partial differential equations. The principal finite elements for differential forms are presented in the periodic table of finite elements \cite{arnold2014}. Along with the shape functions, the table provides degrees of freedom, defined as weighted moments, and together they specify the finite element space on a given mesh. Although these traditional dofs suit the finite element method excellently, for cochain-based methods it is desirable to obtain dofs for $p$-forms through integration on $p$-chains of the mesh. For example, in the case of (lowest order) Whitney forms (i.e.\ the space $\mathcal{P}^-_1\Lambda^p$), the basis $p$-forms are in correspondence with $p$-cochains of the mesh, and hence they can be used as a tool in methods that are based on discrete exterior calculus. With higher order Whitney forms ($\mathcal{P}^-_k\Lambda^p$ for $k>1$) this is no longer the case, and the traditional dofs lack physical interpretation.

Rapetti and Bossavit \cite{rapetti2009} addressed this issue by introducing an approach based on small simplices, which are images of the mesh simplices through homothetic transformations. The idea is to define the shape functions and their dofs using these: to each small $p$-simplex of order $k$ corresponds a Whitney $p$-form of order $k$, and the dofs are obtained through integration over $k$th order small $p$-simplices. Although the approach generalises the lowest order case (in that $k=1$ yields the standard Whitney forms on the initial simplices), the higher order case is not equally simple. In particular, the small simplices do not pave the initial mesh, and the spanning forms corresponding to small simplices are not linearly independent. Despite these downsides, the approach can be reconciled with discrete exterior calculus and has been adopted for use \cite{lohi2019,kettunen2021,lohi2022}.

In this work, we provide an analogous approach for the space $\mathcal{Q}^-_k\Lambda^p$, the (tensor product) finite element space of differential forms on cubical meshes, to which we hereafter refer as “cubical forms” for short. The approach uses small cubes, which are similar to small simplices but defined on cubical meshes. We first give a definition of the small cubes and use them to define cubical forms similarly as higher order Whitney forms are defined using small simplices. The new degrees of freedom resulting from integration over small cubes are considered next: we provide an explicit formula for integrating basis functions and prove that the dofs are unisolvent. Finally, we conclude with the properties of the resulting interpolation operator. Two improvements to the analogous strategy based on small simplices are that the small cubes completely pave the initial mesh and the spanning cubical forms are linearly independent. The approach is hence readily compatible with discrete exterior calculus and enables higher order methods on cubical meshes.

\section{Small cubes and cubical forms}
\label{cubical_forms}

We first define the small cubes and the cubical forms in the unit $n$-cube $\square^n=[0,1]^n$. Cubical meshes are considered in Section~\ref{interpolating}.

\begin{definition}[Small cubes]
Let $\mathcal{J}(n,k-1)$ denote the set of multi-indices $\mathbf{k}=(k_1,\ldots, k_n)$ with $n$ components $k_i\le k-1$. For the unit $n$-cube $\square^n=[0,1]^n$, each multi-index $\mathbf{k}\in\mathcal{J}(n,k-1)$ defines a map $\mathbf{k}_{k-1}:\square^n\to\square^n$ by
\begin{align*}
\mathbf{k}_{k-1}(x_1,\ldots,x_n)=\frac{(k_1+x_1,\ldots,k_n+x_n)}{k}.
\end{align*}
For $k\ge1$, the set of $k$th order small $p$-cubes of $\square^n$ is
\begin{align*}
S_k^p(\square^n)=\{\mathbf{k}_{k-1}(\tau) \mid \mathbf{k}\in\mathcal{J}(n,k-1) \; \textup{and $\tau$ is a $p$-face of $\square^n$}\}.
\end{align*}
\end{definition}

\begin{remark}
Since $\mathcal{J}(n,k-1)\subset\mathcal{J}(n,k)$, the map $\mathbf{k}_{k-1}$ is not defined by the components of $\mathbf{k}$ alone. The subscript specifies the set of multi-indices whose element $\mathbf{k}$ is considered.
\end{remark}

\begin{figure}[htb]
\centering
\includegraphics[width=10cm]{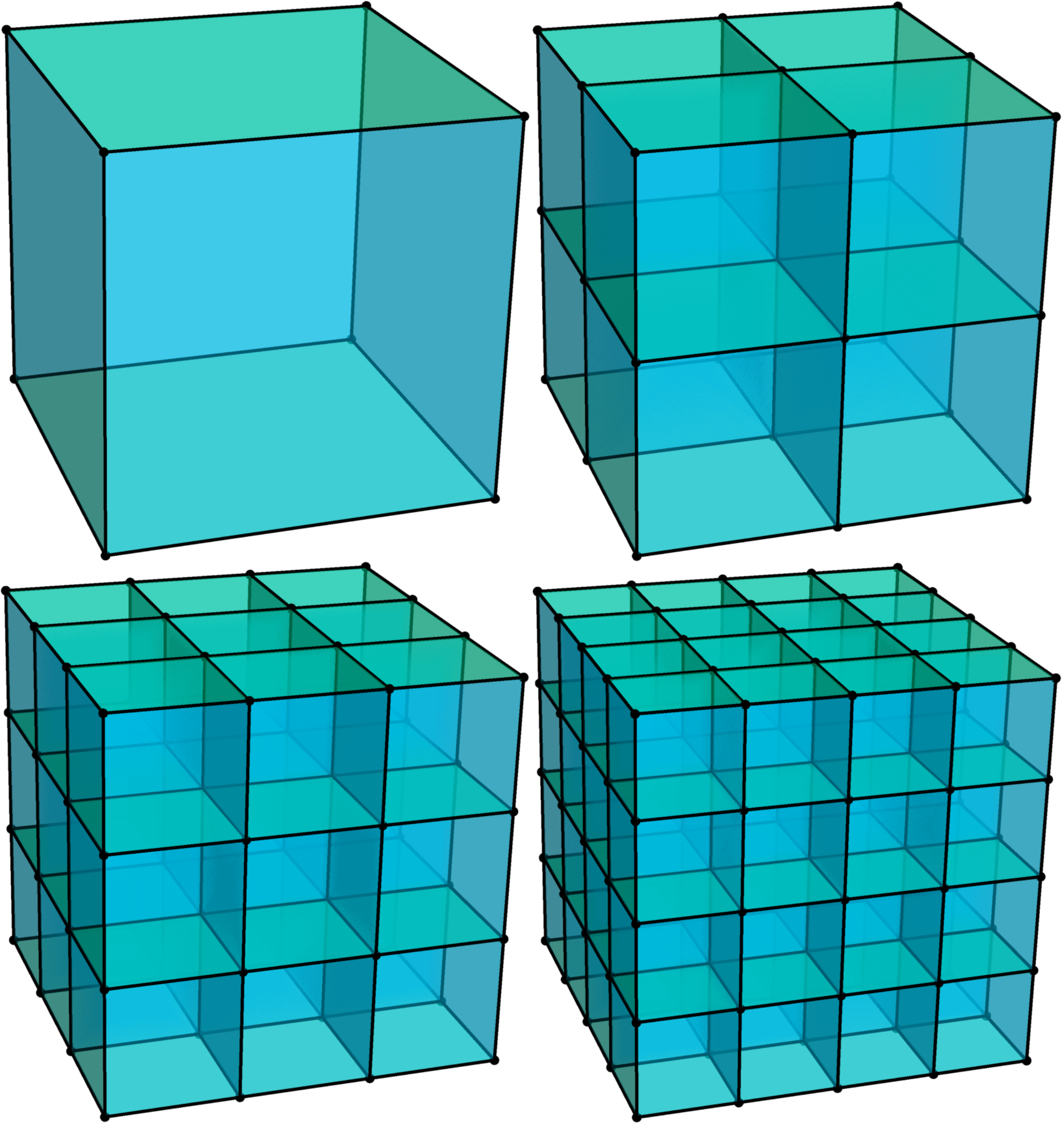}
\caption{Small cubes of orders 1--4 in three dimensions.}
\label{small_cubes_3D}
\end{figure}

\noindent Examples of small cubes are shown in Figure \ref{small_cubes_3D}.

Cubical forms can be seen as counterparts of Whitney forms for cubes. These are the shape functions of the $\mathcal{Q}^-_k\Lambda^p$ family in finite element exterior calculus, and they can be obtained using a tensor product construction \cite{arnold2015b}. We define cubical forms using small cubes similarly as higher order Whitney forms are defined using small simplices. Henceforth, we say that two $p$-cells (or hyperplanes) are parallel if one of them can be moved to the hyperplane of the other by translation.

\begin{definition}[Lowest order cubical forms]
\label{lo_cubical_form_def}
Let $\sigma$ be a $p$-face of $\square^n$. Let $x_{i_1},\ldots,x_{i_p}$ be the coordinates whose plane is parallel to $\sigma$ and $x_{i_{p+1}},\ldots,x_{i_n}$ the other coodinates, whose values $y_{i_{p+1}},\ldots,y_{i_n}$ are either 0 or 1 on $\sigma$. The lowest order cubical form $\wm\sigma$ corresponding to $\sigma$ is
\begin{align*}
\wm\sigma=\bigg(\prod_{j=p+1}^nx_{i_j}^{y_{i_j}}(1-x_{i_j})^{1-y_{i_j}}\bigg)\ed x_{i_1}\wedge\ldots\wedge\ed x_{i_p}.
\end{align*}
\end{definition}

\begin{definition}[Higher order cubical forms]
\label{ho_cubical_form_def}
Let $\mathbf{k}\in\mathcal{J}(n,k-1)$ and $\tau$ be a $p$-face of $\square^n$. The $k$th order cubical $p$-form corresponding to the small cube $\mathbf{k}_{k-1}(\tau)$ is
\begin{align*}
w(\mathbf{k}_{k-1}(\tau))=\bigg(\prod_{i=1}^nx_i^{k_i}(1-x_i)^{k-1-k_i}\bigg)\wm\tau.
\end{align*}
The space of $k$th order cubical $p$-forms is
\begin{align*}
Q_k^p(\square^n)=\textup{span}\{w(\mathbf{k}_{k-1}(\tau))\mid \mathbf{k}\in\mathcal{J}(n,k-1) \; \textup{and $\tau$ is a $p$-face of $\square^n$}\}.
\end{align*}
\end{definition}

Let us first verify that the forms given in Definition~\ref{ho_cubical_form_def} indeed yield the space $\mathcal{Q}^-_k\Lambda^p$.

\begin{proposition}
In the unit $n$-cube $\square^n$, we have $Q_k^p(\square^n)=\mathcal{Q}^-_k\Lambda^p(\square^n)$.
\end{proposition}
\begin{proof}
Recall that $\mathcal{Q}^-_k\Lambda^p(\square^n)$ is spanned by $p$-forms of the form $f\ed x_{i_1}\wedge\ldots\wedge\ed x_{i_p}$, where $f$ is at most $k$th order polynomial in all variables and at most $(k-1)$th order polynomial in the variables $x_{i_1},\ldots,x_{i_p}$. Hence $Q_k^p(\square^n)\subset\mathcal{Q}^-_k\Lambda^p(\square^n)$ follows directly from Definitions~\ref{lo_cubical_form_def} and \ref{ho_cubical_form_def}. It remains to prove $\mathcal{Q}^-_k\Lambda^p(\square^n)\subset Q_k^p(\square^n)$, and for this it is sufficent to show that $Q_k^p(\square^n)$ contains all $p$-forms of the form $x_1^{y_1}\cdot\ldots\cdot x_n^{y_n}\ed x_{i_1}\wedge\ldots\wedge\ed x_{i_p}$, where the $y_i$ are integers such that $0\le y_i\le k$ for all $i$ and $y_i\le k-1$ if $i\in\{i_1,\ldots,i_p\}$.

Let $\omega=x_1^{y_1}\cdot\ldots\cdot x_n^{y_n}\ed x_{i_1}\wedge\ldots\wedge\ed x_{i_p}$ for such integers $y_i$. We choose $z_i=y_i+1$ if $i\in\{i_1,\ldots,i_p\}$, $z_i=y_i$ if $i\notin\{i_1,\ldots,i_p\}$, and write
\begin{align*}
x_i^{y_i}=x_i^{y_i}(x_i+(1-x_i))^{k-z_i}=x_i^{y_i}\sum_{j=0}^{k-z_i}\binom{k-z_i}{j}x_i^j(1-x_i)^{k-z_i-j}.
\end{align*}
Expanding $\omega$ in this way, we get a linear combination of terms of the form
\begin{align*}
\bigg(\prod_{i=1}^nx_i^{a_i}(1-x_i)^{b_i}\bigg)\ed x_{i_1}\wedge\ldots\wedge\ed x_{i_p},
\end{align*}
where $a_i+b_i=k-1$ if $i\in\{i_1,\ldots,i_p\}$ and $a_i+b_i=k$ otherwise. From Definitions~\ref{lo_cubical_form_def} and \ref{ho_cubical_form_def}, we see that such terms are in $Q_k^p(\square^n)$.
\end{proof}

From existing results for $\mathcal{Q}^-_k\Lambda^p$ (see \cite{arnold2015b}), we know that the exterior derivative $\ed$ satisfies $\ed(Q_k^p(\square^n))\subset Q_k^{p+1}(\square^n)$ and the dimension of the space $Q_k^p(\square^n)$ is $\binom{n}{p}k^p(k+1)^{n-p}$. It is easy to see that this is also the number of distinct $k$th order small $p$-cubes of $\square^n$. The spanning forms given in Definition~\ref{ho_cubical_form_def} are hence linearly independent, which is an improvement to the analogous approach based on small simplices and higher order Whitney forms.

\section{New degrees of freedom}

Since $p$-forms can be integrated over small $p$-cubes, we can take the integrals over $k$th order small $p$-cubes as degrees of freedom for $k$th order cubical $p$-forms. Note that each dof can be associated with a specific face of $\square^n$ --- the one that contains the small simplex but has no faces of lower dimension that also contain it. Hence the basic requirement for degrees of freedom is fulfilled: the values of dofs associated with a face only depend on the trace of the differential form on that face.

\subsection{Integrating basis functions over small simplices}

In this subsection we provide a formula for computing the values of the new dofs for basis functions. The following lemmas play a key role.

\begin{lemma}
\label{integral_lemma_1d}
For integers $m,n\ge0$ and for $y,z\in\mathbb{R}$,  
\begin{align*}
\int_0^1(z+x)^n(y+1-x)^mdx=\sum_{i=0}^m\sum_{j=0}^n\binom{m}{i}\binom{n}{j}y^{m-i}z^{n-j}\frac{i!j!}{(i+j+1)!}.
\end{align*}
\end{lemma}
\begin{proof}
\begin{align*}
&\int_0^1(z+x)^n(y+1-x)^mdx=\int_0^1\bigg(\sum_{j=0}^n\binom{n}{j}z^{n-j}\cdot x^j\bigg)\bigg(\sum_{i=0}^m\binom{m}{i}y^{m-i}\cdot(1-x)^i\bigg)dx\\
&=\sum_{i=0}^m\sum_{j=0}^n\binom{m}{i}\binom{n}{j}y^{m-i}z^{n-j}\int_0^1(1-x)^ix^jdx = \sum_{i=0}^m\sum_{j=0}^n\binom{m}{i}\binom{n}{j}y^{m-i}z^{n-j}\frac{i!j!}{(i+j+1)!},
\end{align*}
where we used a well-known integration rule for products of barycentric functions \cite{vermolen2018} in the last step.
\end{proof}

\begin{lemma}
\label{integral_average_lemma}
Let $\tau$ be a $p$-face of $\square^n$. Let $x_{i_1},\ldots,x_{i_p}$ be the coordinates whose plane is parallel to $\tau$ and $x_{i_{p+1}},\ldots,x_{i_n}$ the other coodinates, whose values $y_{i_{p+1}},\ldots,y_{i_n}$ are either 0 or 1 on $\tau$. Let $\mathbf{k}\in\mathcal{J}(n,k)$, $\mathbf{k}'\in\mathcal{J}(n,k')$, and $\upsilon=\mathbf{k}_{k'}'(\tau)$. The average of $\prod_{i=1}^nx_i^{k_i}(1-x_i)^{k-k_i}$ over the small $p$-cube $\upsilon$ is
\begin{align*}
&\frac1{\lvert\upsilon\rvert}\int_{\upsilon}\prod_{i=1}^nx_i^{k_i}(1-x_i)^{k-k_i}=\frac1{(k'+1)^{nk}}\bigg(\prod_{j=p+1}^{n}(k_{i_j}'+y_{i_j})^{k_{i_j}}(k'-k_{i_j}'+1-y_{i_j})^{k-k_{i_j}}\bigg)\\&\cdot\bigg(\prod_{j=1}^p\int_0^1(k_{i_j}'+x)^{k_{i_j}}(k'-k_{i_j}'+1-x)^{k-k_{i_j}}dx\bigg).
\end{align*}
\end{lemma}
\begin{proof}
Recall that $\mathbf{k}_{k'}'$ maps $(x_1,\ldots,x_n)$ to $(k_1'+x_1,\ldots,k_n'+x_n)/(k'+1)$. For $j>p$, $x_{i_j}^{k_{i_j}}(1-x_{i_j})^{k-k_{i_j}}$ has the constant value 
\begin{align*}
\bigg(\frac{k_{i_j}'+y_{i_j}}{k'+1}\bigg)^{k_{i_j}}\bigg(1-\frac{k_{i_j}'+y_{i_j}}{k'+1}\bigg)^{k-k_{i_j}}=\frac1{(k'+1)^k}(k_{i_j}'+y_{i_j})^{k_{i_j}}(k'-k_{i_j}'+1-y_{i_j})^{k-k_{i_j}}
\end{align*}
on $\upsilon$ and hence
\begin{align*}
&\frac1{\lvert\upsilon\rvert}\int_{\upsilon}\prod_{i=1}^nx_i^{k_i}(1-x_i)^{k-k_i}=\frac1{(k'+1)^{(n-p)k}}\cdot\bigg(\prod_{j=p+1}^{n}(k_{i_j}'+y_{i_j})^{k_{i_j}}(k'-k_{i_j}'+1-y_{i_j})^{k-k_{i_j}}\bigg)\\&\cdot\frac1{\lvert\upsilon\rvert}\int_{\upsilon}\prod_{j=1}^px_{i_j}^{k_{i_j}}(1-x_{i_j})^{k-k_{i_j}}.
\end{align*}
Since $1/(k'+1)^p$ is the Jacobian determinant of $\mathbf{k}_{k'}'$ regarded as a map from $\tau$ onto $\upsilon$ and $\frac1{\lvert\upsilon\rvert}=(k'+1)^p$, we can write
\begin{align*}
&\frac1{\lvert\upsilon\rvert}\int_{\upsilon}\prod_{j=1}^px_{i_j}^{k_{i_j}}(1-x_{i_j})^{k-k_{i_j}}=\int_{\tau}\prod_{j=1}^p\bigg(\frac{k_{i_j}'+x_{i_j}}{k'+1}\bigg)^{k_{i_j}}\bigg(1-\frac{k_{i_j}'+x_{i_j}}{k'+1}\bigg)^{k-k_{i_j}}\\&=\frac1{(k'+1)^{pk}}\int_{\tau}\prod_{j=1}^p(k_{i_j}'+x_{i_j})^{k_{i_j}}(k'-k_{i_j}'+1-x_{i_j})^{k-k_{i_j}}.
\end{align*}
The result follows, since the integral above is
\begin{align*}
&\int_{\tau}\prod_{j=1}^p(k_{i_j}'+x_{i_j})^{k_{i_j}}(k'-k_{i_j}'+1-x_{i_j})^{k-k_{i_j}}\\
&=\int_{[0,1]^p}\bigg(\prod_{j=1}^p(k_{i_j}'+x_{i_j})^{k_{i_j}}(k'-k_{i_j}'+1-x_{i_j})^{k-k_{i_j}}\bigg)dx_{i_1}\ldots dx_{i_p}\\
&=\prod_{j=1}^p\int_0^1(k_{i_j}'+x)^{k_{i_j}}(k'-k_{i_j}'+1-x)^{k-k_{i_j}}dx.
\end{align*}
\end{proof}

The integral of any $k$th order spanning $p$-form given in Definition~\ref{ho_cubical_form_def} over any $k$th order small $p$-cube can now be computed by combining Lemmas \ref{integral_lemma_1d} and \ref{integral_average_lemma} with the following proposition.

\begin{proposition}
\label{integral_formula}
Let $\sigma$ be a $p$-face of the unit $n$-cube $\square^n$, and let $\omega$ be a smooth $0$-form. For any small $p$-cube $\upsilon$, we have
\begin{align*}
\int_{\upsilon}\omega\wm\sigma=\bigg(\frac1{\lvert\upsilon\rvert}\int_{\upsilon}\omega\bigg)\langle\wm\sigma(x),\vect(\upsilon)\rangle,
\end{align*}
where $\frac1{\lvert\upsilon\rvert}\int_{\upsilon}\omega$ is the average of $\omega$ over $\upsilon$, $x$ is any point in $\upsilon$, and $\vect(\upsilon)$ is the $p$-vector of $\upsilon$.
\end{proposition}
\begin{proof}
Let $x_{i_1},\ldots,x_{i_p}$ be the coordinates whose plane is parallel to $\sigma$. If $\upsilon$ is not parallel to $\sigma$, then both sides become zero because some of the coordinates is constant and hence $\ed x_{i_1}\wedge\ldots\wedge\ed x_{i_p}$ vanishes on $\upsilon$. But if $\upsilon$ is parallel to $\sigma$, $\wm\sigma$ is constant in $\upsilon$ and hence
\begin{align*}
\int_{\upsilon}\omega\wm\sigma=\int_{\upsilon}\bigg\langle\omega(x)\wm\sigma(x),\frac{\vect(\upsilon)}{\lvert\upsilon\rvert}\bigg\rangle dx=\bigg(\frac1{\lvert\upsilon\rvert}\int_{\upsilon}\omega\bigg)\langle\wm\sigma(x),\vect(\upsilon)\rangle.
\end{align*}
\end{proof}

\subsection{Proof of unisolvence}

Let us next show that these new degrees of freedom are unisolvent. Note that since the number of small $p$-cubes is equal to the number of (linearly independent) spanning $p$-forms, it is sufficient to prove that $\omega\in Q_k^p(\square^n)$ has zero integral over all $k$th order small $p$-cubes only if $\omega=0$. This is shown in Theorem~\ref{unisolvence}, whose proof uses the following two lemmas.

\begin{lemma}
\label{unisolvence_lemma1}
For each $i\in\{1,\ldots,n\}$, let $k_i\ge0$ be an integer and $K_i$ a set of $k_i+1$ distinct real numbers. Suppose that $f:\mathbb{R}^n\to\mathbb{R}$ is a polynomial of order $k_i$ at most in the variable $x_i$, for all $i=1,\ldots,n$. If $f(x)=0$ for all $x\in K_1 \times\ldots\times K_n$, then $f=0$.
\end{lemma}
\begin{proof}
A well-known result for univariate polynomials states that a polynomial of order $k\ge1$ can have at most $k$ roots. Hence the case $n=1$ is clear. Suppose as an induction hypothesis that the statement holds for $n=m-1$, with $m\ge2$, and consider the case $n=m$. If $f(x)=0$ for all $x\in K_1 \times\ldots\times K_m$, then for each $y_j\in K_m$ the function $g_j:\mathbb{R}^{m-1}\to\mathbb{R}$ defined by $g_j(x) = f(x,y_j)$ is zero by the induction hypothesis. Hence for any $(x_1,\ldots,x_{m-1})$, the function $y\mapsto f(x_1,\ldots,x_{m-1},y)$ vanishes in $K_m$ and hence has $k_m+1$ roots. Since it is an univariate polynomial of order $k_m$ at most, it must be zero. Hence the statement holds for $n=m$.
\end{proof}

\begin{lemma}
\label{unisolvence_lemma2}
Suppose that $f:\mathbb{R}^n\to\mathbb{R}$ is a nonzero polynomial. For any $h_1,\ldots,h_n>0$, there exist $\epsilon, M_1,\ldots,M_n>0$ such that $\lvert f\rvert\ge\epsilon$ in $[M_1,M_1+h_1]\times\ldots\times[M_n,M_n+h_n]$.
\end{lemma}
\begin{proof}
We can write
\begin{align*}
f=\sum_{i_1=0}^{k_1}\sum_{i_2=0}^{k_2}\!\ldots\!\sum_{i_n=0}^{k_n}\!a(i_1,i_2,\ldots,i_n)x_1^{i_1}x_2^{i_2}\cdot\ldots\cdot x_n^{i_n}=\sum_{i_1=0}^{k_1}x_1^{i_1}\sum_{i_2=0}^{k_2}x_2^{i_2}\!\ldots\!\sum_{i_n=0}^{k_n}\!a(i_1,i_2,\ldots,i_n)x_n^{i_n},
\end{align*}
where $k_i$ is the order of $f$ in the variable $x_i$ and each coefficient $a(i_1,i_2,\ldots,i_n)$ is constant. For each $j\in\{1,\ldots,n\}$ and $i_1,\ldots,i_{n-j}$ such that $0\le i_l\le k_l$ for all $l\in\{1,\ldots,n-j\}$, let us define a function $g_j[i_1,\ldots,i_{n-j}]:\mathbb{R}^j\to\mathbb{R}$ by
\begin{align*}
g_j[i_1,\ldots,i_{n-j}](x_{n-j+1},\ldots,x_n)=\sum_{i_{n-j+1}=0}^{k_{n-j+1}}x_{n-j+1}^{i_{n-j+1}}\sum_{i_{n-j+2}=0}^{k_{n-j+2}}x_{n-j+2}^{i_{n-j+2}}\ldots\sum_{i_n=0}^{k_n}a(i_1,i_2,\ldots,i_n)x_n^{i_n}.
\end{align*}
In other words, we have
\begin{align*}
&g_1[i_1,\ldots,i_{n-1}](x_n)=\sum_{i_n=0}^{k_n}a(i_1,i_2,\ldots,i_n)x_n^{i_n},\\
&g_j[i_1,\ldots,i_{n-j}](x_{n-j+1},\ldots,x_n)=\sum_{i_{n-j+1}=0}^{k_{n-j+1}}g_{j-1}[i_1,\ldots,i_{n-j+1}](x_{n-j+2},\ldots,x_n)x_{n-j+1}^{i_{n-j+1}},
\end{align*}
and $g_n(x)=f(x)$.

We proceed as follows. At step 1, we can find $\epsilon_n,M_n>0$ such that each $g_1[i_1,\ldots,i_{n-1}]$ is either identically zero or satisfies $\lvert g_1[i_1,\ldots,i_{n-1}](x_n)\rvert\ge\epsilon_n$ for all $x_n\in[M_n,M_n+h_n]$. At step $j$ (for $2\le j\le n$), suppose we have found $\epsilon_{n-j+2},M_{n-j+2},\ldots,M_{n}>0$ such that each $g_{j-1}[i_1,\ldots,i_{n-j+1}]$ is either identically zero or satisfies 
\begin{align*}
\lvert g_{j-1}[i_1,\ldots,i_{n-j+1}](x_{n-j+2},\ldots,x_n)\rvert\ge\epsilon_{n-j+2}
\end{align*}
for all $(x_{n-j+2},\ldots,x_n)\in[M_{n-j+2},M_{n-j+2}+h_{n-j+2}]\times\ldots\times[M_n,M_n+h_n]$. Then we can find $\epsilon_{n-j+1},M_{n-j+1}>0$ such that each $g_j[i_1,\ldots,i_{n-j}]$ is either identically zero or satisfies 
\begin{align*}
\lvert g_j[i_1,\ldots,i_{n-j}](x_{n-j+1},\ldots,x_n)\rvert\ge\epsilon_{n-j+1}
\end{align*}
for all $(x_{n-j+1},\ldots,x_n)\in[M_{n-j+1},M_{n-j+1}+h_{n-j+1}]\times\ldots\times[M_n,M_n+h_n]$. The proof is completed at step $n$, since $g_n=f$, which is nonzero by assumption.
\end{proof}

\begin{theorem}
\label{unisolvence}
Let $\omega\in Q_k^p(\square^n)$. If $\int_{\upsilon}\omega=0$ for all small $p$-cubes $\upsilon\in S_k^p(\square^n)$, then $\omega=0$. 
\end{theorem}
\begin{proof}
Assume $\int_{\upsilon}\omega=0$ for all $\upsilon\in S_k^p(\square^n)$ and write
\begin{align*}
\omega=\sum_{1\le i_1<\ldots<i_p\le n}\omega_{i_1\ldots i_p}\ed x_{i_1}\wedge\ldots\wedge\ed x_{i_p}.
\end{align*}
Note that $\ed x_{i_1}\wedge\ldots\wedge\ed x_{i_p}$ is zero on $\upsilon$ unless $\upsilon$ is parallel to the corresponding coordinate plane. Hence $\int_{\upsilon}\omega=0$ implies $\int_{\upsilon}\omega_{i_1\ldots i_p}\ed x_{i_1}\wedge\ldots\wedge\ed x_{i_p}=0$ for all $1\le i_1<\ldots<i_p\le n$. We show that each coefficient function $\omega_{i_1\ldots i_p}$ is zero.

Let $\tau$ be the $p$-face of $\square^n$ which is parallel to the coordinate plane of $x_{i_1},\ldots,x_{i_p}$ and on which the other coordinates are zero, and let $\mathbf{k}=(0,\ldots,0)\in\mathcal{J}(n,k-1)$. Denote $\upsilon_0=\mathbf{k}_{k-1}(\tau)$ and define a function $f:\mathbb{R}^n\to\mathbb{R}$ by
\begin{align*}
f(u)=\int_{\upsilon_0}\omega_{i_1\ldots i_p}(x+u)\ed x_{i_1}\wedge\ldots\wedge\ed x_{i_p}.
\end{align*}
Observe that since $\omega_{i_1\ldots i_p}$ is at most $k$th order polynomial in all variables and at most $(k-1)$th order polynomial in the variables $x_{i_1},\ldots,x_{i_p}$, the same holds for $f$.

In the small $p$-cube $\upsilon_0$, the coordinates $x_{i_1},\ldots,x_{i_p}$ vary from $0$ to $1/k$ and the other coordinates are zero. The other small $p$-cubes of order $k$ that are parallel to $\upsilon_0$ are obtained from $\upsilon_0$ through translation as follows. Let
\begin{align*}
K_i = \left\{ \begin{array}{cl}
\{0,\frac{1}{k},\frac{2}{k},\ldots,\frac{k-1}{k}\} & \textrm{if $i\in\{i_1,\ldots,i_p\}$},\\
\{0,\frac{1}{k},\frac{2}{k},\ldots,\frac{k-1}{k},1\} & \textrm{if $i\notin\{i_1,\ldots,i_p\}$}.
\end{array} \right.
\end{align*}
Then the small $p$-cubes of order $k$ that are parallel to the coordinate plane of $x_{i_1},\ldots,x_{i_p}$ are precisely the translations of $\upsilon_0$ by vectors $u\in K_1\times\ldots\times K_n$. In particular, we have $f(u)=0$ for all $u\in K_1\times\ldots\times K_n$, and hence $f=0$ by Lemma~\ref{unisolvence_lemma1}.

It remains to show how $f=0$ implies $\omega_{i_1\ldots i_p}=0$. If $\omega_{i_1\ldots i_p}\neq0$, applying Lemma~\ref{unisolvence_lemma2} with $h_1=h_2=\ldots=h_n=1$ yields $\epsilon>0$ and $M_1,\ldots,M_n>0$ such that $\lvert\omega_{i_1\ldots i_p}\rvert\ge\epsilon$ in $[M_1,M_1+1]\times\ldots\times[M_n,M_n+1]$. But $\omega_{i_1\ldots i_p}$ must attain the value 0 somewhere in this set because $f(M_1\ldots,M_n)=0$. This is a contradiction. Hence $\omega_{i_1\ldots i_p}$ must vanish identically, which concludes the proof.
\end{proof}

\section{Interpolating with cubical forms}
\label{interpolating}

Similarly as Whitney forms are used to interpolate cochains on simplicial meshes, cubical forms can be used for interpolating on cubical meshes. We say that a mesh $K$ in $\Omega\subset\mathbb{R}^n$ is cubical if for each $n$-cell $\sigma$ in $K$ there exists an affine bijection $\phi:\square^n\to\sigma$. In other words, we require that $\sigma$ be a parallelotope. (The requirement could be relaxed to accommodate curvilinear meshes, but this would have a negative effect on the approximation properties \cite{arnold2015b}). We denote by $S^p(K)$ the set of $p$-cells and by $C_p^*(K)$ the space of $p$-cochains.

The small cubes of $\sigma\in S^n(K)$ are obtained as the images of the small cubes of $\square^n$ through the map $\phi$, and the corresponding cubical forms in $\sigma$ are defined as the pullbacks through $\phi^{-1}$:
\begin{align*}
w(\phi(\mathbf{k}_{k-1}(\tau)))=(\phi^{-1})^*(w(\mathbf{k}_{k-1}(\tau))).
\end{align*}
When $K$ is a cubical mesh, we define the space of $k$th order cubical $p$-forms as the span of all $k$th order cubical $p$-forms in the cells of $K$. Denote this space by $Q_k^p(K)$. We remark that the space admits a geometric decomposition, in the sense of \cite{arnold2009}, as follows. Let $\mathring{S}_k^p(\sigma^q)$ denote those small $p$-cubes of $\sigma^q\in S^q(K)$ that are not contained in the boundary of $\sigma^q$ and $\mathring{Q}_k^p(\sigma^q)$ those $p$-forms in $Q_k^p(\sigma^q)$ that have zero trace on the boundary of $\sigma^q$. Then
$\mathring{Q}_k^p(\sigma^q)=\textup{span}\{w(\upsilon)\mid \upsilon\in\mathring{S}_k^p(\sigma^q)\}$
and we have the geometric decomposition
\begin{align}
\label{decomposition}
Q_k^p(\sigma^n)=\bigoplus_{\substack{\sigma^q\in S^q(\sigma^n),\\ p\le q\le n}}\mathring{Q}_k^p(\sigma^q), \quad Q_k^p(K)=\bigoplus_{\substack{\sigma^q\in S^q(K),\\ p\le q\le n}}\mathring{Q}_k^p(\sigma^q),
\end{align}
where we have extended elements in $Q_k^p(\sigma^q)$ to elements of $Q_k^p(\sigma^n)$ using a suitable extension operator. (If $\sigma^q$ is a $q$-face of $\sigma^n$, then any small $p$-cube $\upsilon$ of $\sigma^q$ is also a small $p$-cube of $\sigma^n$, so $w(\upsilon)$ extends to $\sigma^n$ by regarding $\upsilon$ as a small $p$-cube of $\sigma^n$.) A dual decomposition can also be obtained by replacing $Q$ and $\mathring{Q}$ in \eqref{decomposition} with $S$ and $\mathring{S}$.

To apply cubical forms with discrete exterior calculus, we refine the cubical mesh $K$ into a finer mesh $K_k$ whose cells are the $k$th order small cubes. Notice that the small cubes pave the initial cubes completely, so there are no holes between them and the refinement $K_k$ is unique (unlike with small simplices). We define the interpolation operator $\im:C_p^*(K_k)\to Q_k^p(K)$ by requiring that
\begin{align}
\label{defining_property}
\int_{\upsilon}\im X=X(\upsilon)
\end{align}
for all $\upsilon\in S^p(K_k)$. The interpolation operator satisfies all expected properties:
\begin{align}
\label{im_property1}
&\drm_k\im X = X \quad \forall X\in C_p^*(K_k),\\
\label{im_property2}
&\im\drm_k\omega=\omega \quad \forall\omega\in Q_k^p(K),\\
\label{im_property3}
&\im\ed X=\ed\im X \quad \forall X\in C_p^*(K_k),
\end{align}
where $\drm_k$ denotes the de Rham map of $K_k$ and $\ed$ denotes both the coboundary operator and the exterior derivative.
\begin{proposition}
The interpolation operator $\im$ is well defined by \eqref{defining_property} and satisfies the properties \eqref{im_property1}--\eqref{im_property3}.
\end{proposition}
\begin{proof}
Theorem~\ref{unisolvence} implies that the restriction of $\drm_k$ to $Q_k^p(K)$ is injective; since the dimensions of $C_p^*(K_k)$ and $Q_k^p(K)$ match, it is bijective, and \eqref{defining_property} defines $\im$ as its inverse. Hence the properties \eqref{im_property1}--\eqref{im_property2} hold. For \eqref{im_property3} we invoke also $\ed(Q_k^p(K))\subset Q_k^{p+1}(K)$: $\im\ed X=\im\ed\drm_k\im X=\im\drm_k\ed\im X=\ed\im X$, where we used \eqref{im_property1}, Stokes' theorem, and the fact that $\ed\im X$ is in $Q_k^{p+1}(K)$.
\end{proof}

\begin{remark}
At this point, we obtain an easy proof for the exact sequence property of cubical forms: if $\Omega$ has trivial homology groups, the spaces $Q_k^p(K)$ constitute an exact sequence with $\ed$. To see this, suppose $\omega\in Q_k^{p+1}(K)$ such that $\ed\omega=0$. Then $\ed\drm_k\omega=\drm_k\ed\omega=0$, and it is a standard result in algebraic topology \cite{whitney1957} that $\drm_k\omega=\ed X$ for some $X\in C_p^*(K_k)$. Hence $\omega=\im\drm_k\omega=\im\ed X=\ed\im X$. It seems that this exact sequence property of cubical forms has not been proven (or even stated) previously in the literature \cite{arnold2013,arnold2015b}.
\end{remark}

The interpolation operator is implemented efficiently using the decomposition~\eqref{decomposition}. To compute the value of $\im X$ in $\sigma^n\in S^n(K)$, we consider basis functions in $\mathring{Q}_k^p(\sigma^q)$ for $q$-faces $\sigma^q$ of $\sigma^n$, with $p\le q\le n$. The coefficients of basis functions in $\mathring{Q}_k^p(\sigma^q)$ only depend on the values of $X$ on those small $p$-cubes that are in $\sigma^q$. Systematic implementation is possible by copying the approach provided in \cite{lohi2022} for higher order Whitney forms and small simplices. With cubical forms the process is only much simpler, since the spanning forms given in Definition~\ref{ho_cubical_form_def} are linearly independent and the refinement $K_k$ has no other cells than small cubes. In addition, now the coefficients of basis functions with $\ed x_{i_1}\wedge\ldots\wedge\ed x_{i_p}$ in them only depend on the values on small cubes that are parallel to the corresponding coordinate plane, which further simplifies the computations.

Besides interpolating cochains, the operator $\im$ can be used to approximate differential forms; the approximation of $\omega$ obtained with cubical forms is $\im\drm_k\omega$. We conclude the paper with a convergence proof for this approximation.
\begin{theorem}
\label{convergence}
Let $\omega$ be a smooth $p$-form in $\Omega$. There exist a constant $C_{\omega,k}$ such that
\begin{align*}
\lvert\im\drm_k\omega(x)-\omega(x)\rvert\le \frac{C_{\omega,k}}{C_{\Theta}^p}h^k\quad\textup{for all $x\in\sigma$ in all $\sigma\in S^n(K)$}
\end{align*}
whenever $h>0$, $C_{\Theta}>0$, and $K$ is a cubical mesh in $\Omega$ such that $\diam(\sigma)\le h$ and $\Theta(\sigma)\ge C_{\Theta}$ for all cells $\sigma$ of $K$.
\end{theorem}

Here $\Theta(\sigma)$ denotes the fullness, which is defined for a $p$-cell $\sigma$ as $\Theta(\sigma)=\lvert\sigma\rvert/\diam(\sigma)^p$. The proof of Theorem~\ref{convergence} is similar to that of Theorem~5.1 in \cite{lohi2021} after some preparations.

\begin{lemma}
\label{parallelotope_lemma}
Let $\sigma$ be an $n$-parallelotope. There exists an $n$-ball $B\subset\sigma$ with diameter $\diam(B)=\Theta(\sigma)\diam(\sigma)$.
\begin{proof}
We may assume that $\sigma=\{\sum_{i=1}^n\mu_iv_i \mid 0\le\mu_i\le1\;\forall i\}$, where $v_1,\ldots,v_n$ are the edge vectors of $\sigma$. Let $\tau$ be any $(n-1)$-face of $\sigma$ and let $h$ denote the distance from the plane of this face to the point $z=\frac12(v_1+\ldots+v_n)$. Since $\lvert\sigma\rvert=2h\lvert\tau\rvert$ and $\lvert\tau\rvert\le\diam(\sigma)^{n-1}$,
\begin{align*}
h=\frac{\lvert{\sigma}\rvert}{2\lvert\tau\rvert}\ge\frac{\lvert{\sigma}\rvert}{2\diam(\sigma)^{n-1}}=\frac12\Theta(\sigma)\diam(\sigma).
\end{align*}
This holds for all $(n-1)$-faces of $\sigma$, and hence the $n$-ball with radius $\frac12\Theta(\sigma)\diam(\sigma)$ centred at $z$ fits in $\sigma$.
\end{proof}
\end{lemma}

Suppose $\sigma\in S^n(K)$ and consider the affine bijection $\phi:\square^n\to\sigma$ from the unit $n$-cube onto $\sigma$. As a consequence of Lemma~\ref{parallelotope_lemma}, we obtain a bound for the norm of $D\phi^{-1}$ as follows. Let $B\subset\sigma$ be an $n$-ball with centre $z$ such that $\diam(B)=\Theta(\sigma)\diam(\sigma)$, and pick $v$ such that $\lvert v\rvert=1$ and $\lvert D\phi^{-1}(z)v\rvert=\max_{\lvert w\rvert=1}\lvert D\phi^{-1}(z)w\rvert$. Since $\phi^{-1}$ is affine, for all $x\in\sigma$
\begin{align}
\label{d_bound}
\begin{split}
\lvert D\phi^{-1}(x)\rvert&=\lvert D\phi^{-1}(z)v\rvert=\frac{\lvert\phi^{-1}(z+\frac12\Theta(\sigma)\diam(\sigma)v)-\phi^{-1}(z-\frac12\Theta(\sigma)\diam(\sigma)v)\rvert}{\Theta(\sigma)\diam(\sigma)}\\&\le\frac{\sqrt{n}}{\Theta(\sigma)\diam(\sigma)},
\end{split}
\end{align}
where we used the fact that $\phi^{-1}$ maps $\sigma$ onto the unit cube, which has diameter $\sqrt{n}$.

\begin{proof}[Proof of Theorem~\ref{convergence}]
We write
\begin{align*}
\omega=\sum_{1\le i_1<\ldots<i_p\le n}\omega_{i_1\ldots i_p}\ed x_{i_1}\wedge\ldots\wedge\ed x_{i_p}
\end{align*}
and, for $y\in\Omega$, denote by $T_{y,i_1\ldots i_p}$ the $(k-1)$th order Taylor polynomial of $\omega_{i_1\ldots i_p}$ at $y$. Since $\omega$ is smooth in $\Omega$, we may find a constant $C_{\omega}$ such that $\lvert\omega_{i_1\ldots i_p}(x)-T_{y,i_1\ldots i_p}(x)\rvert\le C_{\omega}\lvert x-y\rvert^k$ for all $i_1\ldots i_p$ whenever the line segment from $y$ to $x$ is in $\Omega$.

Let $h>0$ and $C_{\Theta}>0$, and suppose $K$ satisfies the assumptions. Fix $\sigma\in K$ and $y\in\sigma$, and denote $g_{i_1\ldots i_p}=\omega_{i_1\ldots i_p}-T_{y,i_1\ldots i_p}$ so that
\begin{align*}
\omega=\sum_{1\le i_1<\ldots<i_p\le n}(T_{y,i_1\ldots i_p}+g_{i_1\ldots i_p})\ed x_{i_1}\wedge\ldots\wedge\ed x_{i_p}, \quad \lvert g_{i_1\ldots i_p}(x)\rvert\le C_{\omega}h^k \; \textup{in $\sigma$}.
\end{align*}
Since $\im\drm_k T_{y,i_1\ldots i_p}\ed x_{i_1}\wedge\ldots\wedge\ed x_{i_p}=T_{y,i_1\ldots i_p}\ed x_{i_1}\wedge\ldots\wedge\ed x_{i_p}$ for all $i_1\ldots i_p$, we have
\begin{align*}
\im\drm_k\omega-\omega=\sum_{1\le i_1<\ldots<i_p\le n}\bigg(\im\drm_k(g_{i_1\ldots i_p}\ed x_{i_1}\wedge\ldots\wedge\ed x_{i_p})-g_{i_1\ldots i_p}\ed x_{i_1}\wedge\ldots\wedge\ed x_{i_p}\bigg).
\end{align*}
In $\sigma$ the interpolant $\im\drm_k(g_{i_1\ldots i_p}\ed x_{i_1}\wedge\ldots\wedge\ed x_{i_p})=\sum_{\upsilon_i\in S_k^p(\sigma)}\alpha_iw(\upsilon_i)$, where each $\alpha_i$ is a linear combination of the integrals of $g_{i_1\ldots i_p}\ed x_{i_1}\wedge\ldots\wedge\ed x_{i_p}$ over small $p$-cubes in $\sigma$. The coefficients in this linear combination are constant and independent of $\sigma$, so we may find a constant $C_{\alpha}$, depending only on $n$, $p$, and $k$, such that for all the coefficients 
\begin{align*}
\lvert\alpha_i\rvert\le C_{\alpha}\max_{\upsilon_j\in S_k^p(\sigma)}\lvert \int_{\upsilon_j}g_{i_1\ldots i_p}\ed x_{i_1}\wedge\ldots\wedge\ed x_{i_p}\rvert\le C_{\alpha}C_{\omega}h^k\cdot\max_{\upsilon_j\in S_k^p(\sigma)}\lvert\upsilon_j\rvert\le C_{\alpha}C_{\omega}h^k\diam(\sigma)^p.
\end{align*}

In the unit cube, we clearly have $\lvert w(\upsilon)(x)\rvert\le1$ for all $\upsilon\in S_k^p(\square^n)$ and $x\in\square^n$. Applying the pullback inequality $\lvert f^*\omega(x)\rvert\le\lvert Df(x)\rvert^p\cdot\lvert\omega(f(x))\rvert$ \cite[II, 4.12]{whitney1957} to the inverse of the affine bijection $\phi:\square^n\to\sigma$ and using \eqref{d_bound}, we get
\begin{align*}
\lvert w(\upsilon_i)(x)\rvert\le\lvert D\phi^{-1}(x)\rvert^p\le\frac{\sqrt{n}^p}{\Theta(\sigma)^p\diam(\sigma)^p}, \quad \lvert\alpha_iw(\upsilon_i)(x)\rvert\le\frac{C_{\alpha}C_{\omega}\sqrt{n}^ph^k}{\Theta(\sigma)^p}
\end{align*}
for all $\upsilon_i\in S_k^p(\sigma)$ and $x\in\sigma$. Hence
\begin{align*}
\lvert\im\drm_k(g_{i_1\ldots i_p}\ed x_{i_1}\wedge\ldots\wedge\ed x_{i_p})(x)\rvert\le\binom{n}{p}k^p(k+1)^{n-p}\frac{C_{\alpha}C_{\omega}\sqrt{n}^ph^k}{\Theta(\sigma)^p}
\end{align*}
and we may choose
\begin{align*}
C_{\omega,k}=\binom{n}{p}\bigg(\binom{n}{p}k^p(k+1)^{n-p}C_{\alpha}C_{\omega}\sqrt{n}^p+C_{\omega}\bigg)
\end{align*}
to obtain
\begin{align*}
&\lvert\im\drm_k\omega(x)-\omega(x)\rvert\le\!\!\!\sum_{1\le i_1<\ldots<i_p\le n}\!\!\!\!\!\lvert\im\drm_kg_{i_1\ldots i_p}\ed x_{i_1}\wedge\ldots\wedge\ed x_{i_p}(x)\rvert+\lvert g_{i_1\ldots i_p}\ed x_{i_1}\wedge\ldots\wedge\ed x_{i_p}(x)\rvert\\
&\le\binom{n}{p}\bigg(\binom{n}{p}k^p(k+1)^{n-p}\frac{C_{\alpha}C_{\omega}\sqrt{n}^ph^k}{\Theta(\sigma)^p}+C_{\omega}h^k\bigg)\le\frac{C_{\omega,k}}{\Theta(\sigma)^p}h^k
\end{align*}
for all $x\in\sigma$ in all $\sigma\in S^n(K)$ whenever $K$ satisfies the assumptions.
\end{proof}

\bibliographystyle{plain}
\bibliography{references}

\begin{thebibliography}{10}

\bibitem{arnold2014}
Douglas Arnold and Anders Logg.
\newblock Periodic table of the finite elements.
\newblock {\em SIAM News}, 47(9), 2014.

\bibitem{arnold2013}
Douglas~N. Arnold.
\newblock Spaces of finite element differential forms.
\newblock In {\em Analysis and Numerics of Partial Differential Equations},
  volume~4 of {\em Springer INdAM Series}, pages 117--140. Springer, 2013.

\bibitem{arnold2015b}
Douglas~N. Arnold, Daniele Boffi, and Francesca Bonizzoni.
\newblock Finite element differential forms on curvilinear cubic meshes and
  their approximation properties.
\newblock {\em Numerische Mathematik}, 129(1):1--20, 2015.

\bibitem{arnold2006}
Douglas~N. Arnold, Richard~S. Falk, and Ragnar Winther.
\newblock Finite element exterior calculus, homological techniques, and
  applications.
\newblock {\em Acta Numerica}, 15:1--155, 2006.

\bibitem{arnold2009}
Douglas~N. Arnold, Richard~S. Falk, and Ragnar Winther.
\newblock Geometric decompositions and local bases for spaces of finite element
  differential forms.
\newblock {\em Computer Methods in Applied Mechanics and Engineering},
  198(21-26):1660--1672, 2009.

\bibitem{kettunen2021}
Lauri Kettunen, Jonni Lohi, Jukka R{\"a}bin{\"a}, Sanna M{\"o}nk{\"o}l{\"a},
  and Tuomo Rossi.
\newblock Generalized finite difference schemes with higher order {Whitney}
  forms.
\newblock {\em ESAIM: Mathematical Modelling and Numerical Analysis}, 55(4),
  2021.

\bibitem{lohi2019}
Jonni Lohi.
\newblock Discrete exterior calculus and higher order {Whitney} forms.
\newblock Master's thesis, University of Jyv{\"a}skyl{\"a}, 2019.

\bibitem{lohi2022}
Jonni Lohi.
\newblock Systematic implementation of higher order {Whitney} forms in methods
  based on discrete exterior calculus.
\newblock {\em Numerical Algorithms}, 91(3):1261--1285, 2022.

\bibitem{lohi2021}
Jonni Lohi and Lauri Kettunen.
\newblock Whitney forms and their extensions.
\newblock {\em Journal of Computational and Applied Mathematics}, 393:113520,
  2021.

\bibitem{rapetti2009}
Francesca Rapetti and Alain Bossavit.
\newblock Whitney forms of higher degree.
\newblock {\em SIAM Journal on Numerical Analysis}, 47(3):2369--2386, 2009.

\bibitem{vermolen2018}
F.~J. Vermolen and A.~Segal.
\newblock On an integration rule for products of barycentric coordinates over
  simplexes in $\mathbb{R}^n$.
\newblock {\em Journal of Computational and Applied Mathematics}, 330:289--294,
  2018.

\bibitem{whitney1957}
Hassler Whitney.
\newblock {\em Geometric Integration Theory}.
\newblock Princeton University Press, 1957.

\end{thebibliography}

\end{document}